\documentclass[12pt, leqno]{amsart}
\usepackage{amsmath, amsthm, amssymb, amscd}
\usepackage{fancyhdr}
\usepackage{mathdots}
\usepackage{setspace}
\usepackage[small,nohug,heads=vee]{diagrams}
\usepackage{young}
\usepackage{epic}
\diagramstyle[labelstyle=\scriptstyle]
\usepackage[all]{xy}
\usepackage{array}
\usepackage{longtable}

\theoremstyle{plain}

\newtheorem{thm}{Theorem}

\newtheorem*{thmnonum}{Theorem}
\newtheorem*{propnonum}{Proposition}
\newtheorem{prop}[thm]{Proposition}
\newtheorem{lem}[thm]{Lemma}
\newtheorem*{quest}{Question}
\newtheorem{corl}[thm]{Corollary}

\theoremstyle{definition}\newtheorem{mydef}[thm]{Definition}
\theoremstyle{remark}\newtheorem{rmk}[thm]{Remark}
\newenvironment{prf}{\textbf{Proof:}}{\rule{1.6ex}{1.6ex}\\}

\date{\today}

\newcommand{\Z}{\mathbb{Z}}
\newcommand{\C}{\mathbb{C}}

\newcommand{\g}{\mathfrak{g}}

\newcommand{\tf}{\mathfrak{t}}

\newcommand{\Lc}{\mathcal{L}}
\newcommand{\Rc}{\mathcal{R}}
\newcommand{\m}{\setminus}
\newcommand{\Lam}{\Lambda}
\newcommand{\lam}{\lambda}

\newcommand{\Del}{\Delta}

\newcommand{\Hom}{\mbox{Hom}}

\newcommand{\bfr}{\mathfrak{b}}

\begin{document}

\title[Descent of line bundles]{Some conditions for descent of line bundles to GIT quotients $(G/B \times G/B \times G/B)//G$}

\author{Nathaniel Bushek}

\address{Department of Mathematics and Statistics\\University of Alaksa-Anchorage\\Anchorage, AK 99508}

\email{nbushek$@$alaska.edu}

\begin{abstract} 
We consider the descent of line bundles to GIT quotients of products of flag varieties. Let $G$ be a simple, connected, algebraic group over $\mathbb{C}$. We fix a Borel subgroup $B$ and consider the diagonal action of $G$ on the projective variety $X = G/B \times G/B \times G/B$. For any triple $(\lam, \mu, \nu)$ of dominant regular characters there is a $G$-equivariant line bundle $\mathcal{L}$ on $X$. Then, $\mathcal{L}$ is said to descend to the GIT quotient $\pi:[X(\mathcal{L})]^{ss} \rightarrow X(\mathcal{L})//G$ if there exists a line bundle $\hat{\mathcal{L}}$ on $X(\mathcal{L})//G$ such that $\mathcal{L}\mid_{[X(\mathcal{L})]^{ss}} \cong \pi^*\hat{\mathcal{L}}$.

Let $Q$ be the root lattice, $\Lambda$ the weight lattice, and $d$ the least common multiple of the coefficients of the highest root $\theta$ of the Lie algebra $\mathfrak{g}$ of $G$ written in terms of simple roots. We show that $\mathcal{L}$ descends if $\lam, \mu, \nu \in d\Lambda$ and $\lam + \mu + \nu \in \Gamma$, where $\Gamma$ is a fixed sublattice of $Q$ depending only on the type of $\mathfrak{g}$. Moreover, $\Lc$ never descends if $\lam + \mu + \nu \notin Q$. 

\end{abstract}

\maketitle

\section{Introduction}\label{intro}
Let $G$ be a simple, connected, complex linear algebraic group. Let $B$ be a fixed Borel subgroup and $T \subset B$ a fixed maximal torus. Let $\g$ be the Lie algebra of $G$, and $\tf$, $\bfr$, the Lie algebras of $T$ and $B$. Let $\Lam$ be the weight lattices of $\g$, $X(T)$ the character group of $T$, $X(T)^{++}$ the dominant regular characters, and $Q$ the root lattice of $G$.

Then $X := (G/B \times G/B \times G/B)$ is a projective variety with a natural action of $G$ given by the diagonal of left multiplication. Let $\Lc$ be an ample line bundle on $X$. Ample line bundles on $X$ correspond to triples of dominant regular characters on $T$. Let $\Lc = \Lc(\lam, \mu, \nu)$ be the ample line bundle associated to $(\lam, \mu, \nu)$, and $X//G$ be the geometric invariant theory quotient of $X$ by $G$ determined by $\Lc$. We consider the following question.

\begin{quest}

What conditions can be placed on a triple of dominant regular characters $(\lam, \mu, \nu)$ to know that the corresponding line bundle $\Lc$ on $X$ will or will not descend to the GIT quotient $X//G$?
\end{quest}

The primary model for our work comes from \cite{kumardescent}. Here, Kumar considers the descent question for the GIT quotient $(G/P)//T$, where $B \subset P$ is a parabolic subgroup and the torus action is by left multiplication. In this case, a line bundle $\Lc(\lam)$ descends if and only $\lam$ belongs to a explicitly given sublattice $\Gamma$ of $Q$ (Table \ref{table:Gammalattice}). The lattice $\Gamma$ continues to play an important role in our case. Let $d$ be the least common multiple of the coefficients of the highest root $\theta$ of $\g$ when expressed in terms of the simple roots (Table \ref{table:theta}). The following theorem is our main sufficient condition which appears in Theorem \ref{56}.

\begin{thmnonum}

Given $\lam, \mu, \nu  \in X(T)^{++}$, if $\lam, \mu, \nu \in  d \Lam$ and $\lam + \mu + \nu \in \Gamma$, then $\Lc(\lam, \mu, \nu)$ descends to $X//G$. 
\end{thmnonum}

To prove the theorem, we utilize Kempf's `descent lemma' (see \cite{drezet}, Theorem 2.3), which reduces the question of descent to the action of stabilizers on fibers of the line bundle. Our approach is similar to that taken in \cite{kumardescent}. However, the case of $X//G$ quickly becomes complicated. For example, in \cite{kumardescent}, stabilizers have elegant and useful descriptions. In our case, stabilizers are given as intersections of conjugates of $B$. These are generally difficult objects. A simple yet important step in our approach is reducing the application of Kempf's lemma to points $x \in X$ with a reductive stabilizer. 

Assumptions are artificially imposed that keep Theorem \ref{56} from being optimal. Section \ref{counterexample} gives a counter example that the conditions in Theorem \ref{56} are not necessary. However, it remains unclear if Theorem \ref{56} provides the best condition that is independent of where $\Lc$ sits in the $G$-ample cone, denoted $C^G(X)$, as defined in \cite{dh} (or, equivalently where $(\lam, \mu, \nu)$ sits in the saturated tensor cone as given in \cite{belkalekum}). Recall, the cone of $G$-ample line bundles on $X$ that admit a non-empty semistable locus is partitioned into GIT classes; within each GIT class, the semistable and stable loci remain constant. Certainly, some GIT classes have sharper descent conditions. It remains unclear if the union of optimal descent conditions on each GIT class will be an improvement of Theorem \ref{56}.

We also obtain the following necessary condition in Proposition \ref{100}. In type $A$, Theorem \ref{56} and Proposition \ref{100} agree, but not in general. 
\begin{propnonum}

If $\lam + \mu + \nu \notin Q$, then $\Lc(\lam, \mu, \nu)$ does not descend. 
\end{propnonum}

While we consider the study of descending line bundles to GIT quotients to be of independent interest (see e.g., \cite{kumarzero}, \cite{pauly}, \cite{teleman}), there is also a major motivation for considering descent in our context. When $\Lc =\Lc(\lam, \mu, \nu)$ descends to a line bundle $\hat{\Lc}$ on $X//G$, we have $H^0(X//G, \hat{\Lc}) \cong H^0(X, \Lc)^G$ (see e.g., \cite{teleman} Theorem 3.2.a). Using the Borel-Weil theorem, the dimension (over $\C$) of the right hand side is $\dim [V(\lam) \otimes V(\mu) \otimes V(\nu)]^G$, which is exactly the multiplicity of the irreducible representation $V(\lam)^*$ inside $V(\mu)\otimes V(\nu)$. On the other hand, due to the vanishing of higher cohomology, $H^0(X//G, \hat{\Lc})$ is the Euler-Poincar\'{e} characteristic of $\hat{\Lc}$ over $X//G$. By the Riemann-Roch Theorem for singular varieties, one can then see that this value is a piecewise polynomial in terms of $\lam, \mu, \nu$. 

In this sense, we would say that the tensor product multiplicity function is \textit{piecewise polynomial}. The general approach to proving piecewise polynomiality when descent occurs is understood and follows as in \cite{kumarzero}, where descent to $(G/P)//T$ was used to prove piecewise polynomiality of the dimension of the zero weight space of a representation. Two obstacles remain in our setting. First, the sectors of polynomiality are not yet fully known. This is the same as describing the GIT classes of maximal dimension in $C^G(X)$. It is worthwhile to note that since $C^G(X)$ is the saturated tensor cone, the boundary of $C^G(X)$ has been fully described by a non-redundant list of inequalities (see \cite{belkalekum}, \cite{ress}). Second, except for type A, our descent conditions are not yet sharp. So, piecewise polynomiality can only be stated on an appropriate sublattice of $C^G(X)$. Ongoing work is aimed at resolving these two obstacles. 

Finally, we note that the piecewise polynomiality of the tensor product multiplicity function is indeed already known. This fact follows separately from the work of Berenstein and Zelevinsky in \cite{bz} and Meinrenken and Sjamaar in \cite{meinsja}. However, Berenstein and Zelevinsky prove that the tensor product multiplicity is given by the number of integral points of some convex polytope, from which piecewise polynomiality follows. Meinrenken and Sjamaar use sympletic geometry. Since our approach aims to give an explicit construction of the polynomial using algebraic-geometry, along with an explicit description of the sectors of polynomiality, we feel that this will be a worthwhile contribution to the existing literature.  

\vskip2ex
 \noindent
{\bf Acknowledgements:} The author would like to express great appreciation to S. Kumar for both suggesting this problem and his helpful guidance at many junctures throughout this work.  

\section{Notation and Preliminaries}\label{notations}

We expand upon the notation presented in the introduction. The Weyl group of $G$ is $W := N_G(T)/T$, where $N_G(T)$ is the normalizer of $T$ in $G$. The set of roots will be denoted $R$, with positive roots $R^+$ (resp. negative roots $R^-$) determined by $B$, simple roots $\Del = \{\alpha_1, \dots, \alpha_\ell\}$, and root lattice $Q := \Z \Del$. Let $\Lam$ be the weight lattice of $\g$ spanned by the fundamental weights $\{\varpi_1, \cdots, \varpi_\ell \}$. If a weight $\lam \in \Lam$ is a positive (resp. non-negative) linear combination of fundamental weights, $\lam$ is a dominant regular (resp. dominant) weight, the set of which is denoted by $\Lam^{++}$ (resp. $\Lam^{+}$). Last, let $X(T)$ be the character group of $T$, i.e., the group of all algebraic group homomorphisms $T \rightarrow \C^*$, and so the set of dominant regular characters is $X(T)^{++}:= X(T) \cap \Lam^{++}$. Recall that $Q \subset X(T) \subset \Lam$. When we are considering a weight $\lam \in \Lam$ as a character, we will write $e^\lam$ to emphasize $e^\lam$ as a homomorphsim $T \rightarrow \C^*$.

We set $U$ to be the unipotent radical of $B$ and $U^-$ to be the unipotent radical of the opposite Borel subgroup $B^-$. Then, any character $e^\chi \in X(T)$ extends to a character on $B$ by setting $e^\chi|_U \equiv 1$ (recall the decomposition $B = TU$). On $G/B$, we have the $G$-equivariant line bundle, $\Lc(\lam) = G \times_B \C_{\lam^{-1}}$, associated to the principle $B$-bundle $G \mapsto G/B$, via the one dimensional $B$-representation $\C_{\lam^{-1}}$. An element in the fiber of $\Lc(\lam)$ over $gB \in G/B$ is the $B$-class $[g,z]$, where $z \in \C^*$ and $[g,z] = [gb^{-1}, bz] = e^{-\lam}(b) [gb^{-1}, z]$ for all $b \in B$. Now, given $\lam, \mu, \nu \in X(T)^{++}$,  taking exterior tensor product gives the $G$-equivariant ample line bundle 
\[
\begin{array}{ll}
 \Lc(\lam, \mu, \nu) & :=  \Lc(\lam) \boxtimes \Lc(\mu) \boxtimes \Lc(\nu) \\ \\
                             & \; =  p_1^*\big[\Lc(\lam)\big] \otimes p_2^*\big[\Lc(\mu)\big] \otimes p_3^*\big[\Lc(\nu)\big] 
\end{array}
\]
on $X$. Here, $p_i: X \rightarrow G/B$ is projection onto the $i$th factor of the product. $G$ acts diagonally on $\Lc(\lam, \mu, \nu)$. When the context is clear, we will assume the depedence on characters and simply write $\Lc$ for $ \Lc(\lam, \mu, \nu)$.

Recall that for a projective $G$-variety with an ample $G$-linearized line bundle $\Lc$ over $X$, the semistable locus is defined as
\[
 X^{ss}(\Lc) := \{ x\in X\;|\; \exists \;\sigma \in H^0(X, \Lc^{\otimes N})^G\;\;\mbox{such that} \;\;\sigma(x) \neq 0\},
\]
for some $N>0$. This clearly depends on the choice of $G$-linearized line bundle, although any two line bundles in the same GIT equaivalence class give the same semistable locus (see e.g., \cite{dh}). Geometric invariant theory provides the existence of a projective variety $X//G$, called the GIT quotient, and a surjective morphism $\pi : X^{ss} \rightarrow X//G$ that is a good quotient (see e.g. \cite{git}). Now, we come to a central definition.

\begin{mydef}
We say that a line bundle $\Lc$ on $X$ \textit{descends to a line bundle on $X//G$} if there exists a line bundle $\hat{\Lc}$ on $X//G$ such that $\pi^*(\hat{\Lc}) \cong \Lc|_{X^{ss}}$, where the isomorphism is $G$-equivariant.
\end{mydef}

If $X^{ss} = \emptyset$, then descent occurs trivially. Therefore, for the remainder of this paper we only consider $(\lam, \mu, \nu) \in  (X(T)^{++})^3$ such that $X^{ss} (\Lc) \neq \emptyset$. From the perspective of geometric invariant theory, the set of such characters is contained in the $G$-ample cone, $C^G(X)$, as defined in \cite{dh}. From the perspective of representation theory, the set of such characters is the saturated tensor semigroup as studied in \cite{belkalekum}.

\section{Reduction to Stabilizers in $T$}

For this section we consider the more general setting of $Y$ any projective $G$-variety with $\Lc$ any $G$-equivariant ample line bundle on $Y$. We maintain our assumptions on $G$ as set earlier. Recall the following `descent' lemma of Kempf (\cite{drezet}, Theorem 2.3) partially adapted to our setting.

\begin{lem}\label{kempfsdescentlemma} $\Lc$ descends to $Y//G$ if and only if for any $y \in Y^{ss}$, the isotropy subgroup $G_y$ acts trivially on the fiber $\Lc_y$. In fact, for the `if' part, it suffices to assume that $G_y$ acts trivially for only those $y \in Y^{ss}$ such that the orbit $G \cdot y$ is closed in $Y^{ss}$. 
\end{lem}

We begin with an important reduction. 

\begin{lem}\label{reductiveaction} If $G_y$ acts trivially on $\Lc_y$ for every $y \in Y^{ss}$ such that $G_y$ is reductive, then $\Lc$ descends. 
\end{lem}

\begin{prf} From Lemma \ref{kempfsdescentlemma}, it suffices to show that $G_y$ acts trivially on $\Lc_y$ for all $y \in Y^{ss}$ such that $G\cdot y $ is closed in $Y^{ss}$. By \cite{dh}, Lemma 3.3.12, if $G\cdot y$ is closed in $Y^{ss}$, then, $G_y$ is reductive (this is essentially Matsushima's Theorem). \end{prf}

\begin{prop}\label{54} Let $H$ be any reductive subgroup of $B$. Then, $H$ is contained in some torus, and $bHb^{-1} \subset T$ for some $b \in B$. 
\end{prop}

\begin{prf} Recall the projection $p_T: B  \rightarrow T$ sending $tu $ to $t$, which is a homomorphism of algebraic groups. Compose $p_T$ with the inclusion $H \subset B$ to get a homomorphism $\psi: H \rightarrow T$. Since $U$ is normal in $B$, $H \cap U$ is a normal, connected, subgroup of $H$. Then, $H$ being reductive, $H \cap U = \{e\}$. So, $\psi$ is injective, and $H$ is isomorphic to a closed subgroup of $T$. Hence, $H$ is contained in a torus, and $bHb^{-1} \subset T$ for some $b \in B$ (see e.g., \cite{LAG}, Corollary 6.3.6 and Theorem 6.3.5).\end{prf}

The following fact is well known and is straightforward to see. 

\begin{lem}\label{orbitlem} If $G_y$ acts trivially on $\Lc_y$, then for any $g \in G$, $G_{gy}$ acts trivially on $\Lc_{gy}$.
\end{lem}

For the remainder of this paper, for $y \in Y$, we set 
\[
 T_y := G_y \cap T.
\]
This subgroup is central to our study, and the following corollary indicates why. We place an additional assumption on $G$ and $Y$ in the following corollary. However, this assumption is easily satisfied in our eventual context.

\begin{corl}\label{21} Assume $G$ acts on $Y$ such that for all $y \in Y$ there is some $g \in G$ with $G_y \subset gBg^{-1}$. If $G_y$ acts trivially on $\Lc_y$ for all $y \in Y^{ss}$ such that $G_y = T_y$, then $\Lc$ descends to $Y//G$.  
\end{corl}

\begin{prf} By Lemmas \ref{reductiveaction} and \ref{orbitlem} it suffices to show that for any $y' \in Y$ with $G_{y'}$ reductive, there is some $g \in G$ such that $G_{gy'} = T_{gy'}$. But, by Proposition \ref{54}, $bG_{y'}b^{-1} = G_{by'} \subset T$ for some $b \in B$. \end{prf}

\section{Structure of Stabilizers in $T$}

Now, we return to the setting of $X = (G/B \times G/B \times G/B)$ from sections \ref{intro} and \ref{notations}. By Lemma \ref{orbitlem}, within each $G$-orbit we can consider the action of the stabilizer at preferred points. Therefore, we restrict our attention to points of the form $x = (B, g_1B, g_2B) \in X$. It is then easy to see that
\begin{equation*}\label{isotropy}
G_x = B \cap g_1Bg_1^{-1} \cap g_2Bg_2^{-1}.
\end{equation*}
In this case, we clearly have $G_x \subset B$. Reductivity of $G_x$ is an obrit-wise condition since conjugates of a reductive group are reductive. Therefore, we can apply Corollary \ref{21} and restrict our attention to points of the form $x = (B, g_1B, g_2B) \in X$ with $G_x = T_x$.

Recall that given any $\beta \in R$, there is an isomorphism $u_\beta: \C \rightarrow U_\beta$ such that $du_\beta(\C) = \g_\beta$, where $\g_\beta$ is the $\beta$ root subspace of $\g$ (see e.g, \cite{LAG}). Given any ordering of negative roots $(\beta_1, \dots, \beta_m)$, we get an isomorphism of varieties $\phi: \C^m \rightarrow U^-$ given by $\phi(x_1, \dots, x_m) = u_{\beta_1}(x_1) \cdots u_{\beta_m}(x_m)$. We fix such an isomorphism $\phi$ for the remainder of this paper. 

For any $w \in W$, let $R(w) = \{ \beta \in R^+\;|\; w \beta <0\}$ be the inversion set of $w$. Recall the following form of the Bruhat decomposition
\[
G = \coprod_{w \in W} U_{w^{-1}}\dot{w} B = \coprod_{w \in W} \dot{w}U_{w^{-1}}^- B, 
\]
where 
\[
U_{w^{-1}} = \prod_{\alpha \in R(w^{-1})} U_{\alpha} \subset U,
\]
\[
U_{w^{-1}}^- = \dot{w}^{-1}U_{w^{-1}}\dot{w}= \prod_{\alpha \in R(w^{-1})} U_{w^{-1}\alpha} \subset U^-,
\] 
and $\dot{w} \in N_G(T)$ denotes any lift of $w \in W$. Thus, any $x = (B, g_1B, g_2B) \in X$ can be written uniquely as 
\begin{equation}\label{niceform}
 (B, g_1B, g_2B) =(B, \dot{w}_1u_1 B, \dot{w}_2 u_2 B)
\end{equation}
with $w_1, w_2 \in W$ and $u_i \in U_{w_i^{-1}}^-$ for $i = 1, 2$. We will regularly consider points $x \in X$ in this form. 

\begin{lem}\label{30} Let $x \in X$ be as in equation (\ref{niceform}). If $t \in T$ satisfies
\begin{equation}\label{toricform}
t =  \dot{w}_1u_1t_1 v_1 u_1^{-1} \dot{w}_1^{-1} = \dot{w}_2 u_2 t_2 v_2 u_2^{-1} \dot{w}_2^{-1}
\end{equation}
for some $t_1, t_2 \in T$ and $ v_1, v_2 \in U$, then $t = \dot{w}_1t_1  \dot{w}_1^{-1} = \dot{w}_2t_2  \dot{w}_2^{-1}$ and $v_1 = v_2 = e$.
\end{lem}

\begin{prf} We prove only $t = \dot{w}_1t_1  \dot{w}_1^{-1}$ as the second proof is identical. By the hypothesis
\begin{equation*}
t_1 v_1 = u_1^{-1}(\dot{w}_1^{-1}t\dot{w}_1)u_1 = (\dot{w}_1^{-1}t\dot{w}_1)(\dot{w}_1^{-1}t\dot{w}_1)^{-1}u_1^{-1}(\dot{w}_1^{-1}t\dot{w}_1)u_1 .
\end{equation*}
Now, $\dot{w}_1^{-1}t\dot{w}_1 \in T$ and $(\dot{w}_1^{-1}t\dot{w}_1)^{-1}u_1^{-1}(\dot{w}_1^{-1}t\dot{w}_1) \in U^{-}$ since $T$ normalizes $U^-$. Thus, the right hand side of the equation is in $B^-$ and the left hand side is in $B$. Because of the decompositions $B= TU$, $B^- = TU^-$, and $B \cap B^- = T$, we conclude $t_1 = (\dot{w}_1^{-1}t\dot{w}_1)$ and $v_1 = e$. \end{prf}

\begin{rmk}\label{intersectremark} In Lemma \ref{30}, the conclusion that $v_1 = e$ for any $t_1$, $v_1$, $w_1$, and $u_1$ satisfying equation (\ref{toricform}) implies that $T \cap u_1 B u_1^{-1} = T \cap  u_1 T u_1^{-1}$.
\end{rmk}

Now, we consider the structure of $T_x$. First, given $u \in U^-$, we define the \textit{roots supporting $u$} to be the set
\begin{equation}
 \Rc(u) := \{\beta_i \in R^- \;|\; \pi_i(\phi^{-1}(u)) \neq 0\;\}
\end{equation}
where $\phi: \C^m \rightarrow U^-$ is the fixed isomorphism of varieties and $\pi_i: \C^m \rightarrow \C$ is the $i$th projection. In other words, $\Rc(u)$ is the set of roots $\beta \in R^-$ for which the root subgroup $u_\beta$ is a non-identity factor of $u$ as expressed in the image of $\phi$. Then, for  $x = (B, \dot{w}_1u_1 B, \dot{w}_2 u_2 B) \in X$ as in equation (\ref{niceform}), we define the following set of roots
\begin{equation}\label{rootsofx}
 R_x := w_1 \Rc(u_1) \cup w_2 \Rc(u_2).
\end{equation}

\begin{rmk}
Since we are taking $u_1 \in U_{w_1^{-1}}^- \subset U^-$ we see that $\Rc(u_1) \subset w_1^{-1}R(w_1^{-1})$, hence $w_1 \Rc(u_1) \subset R(w_1^{-1})$, and similarly for $u_2$. Also, when $u_1$ is a generic element of $U_{w_1^{-1}}^-$, then $w_1 \Rc(u_1) = R(w_1^{-1})$, and similarly for $u_2, w_2$. Thus, we conclude $R_x \subset R(w_1^{-1}) \cup R(w_2^{-1})$ for any $x$ in the form of equation (\ref{niceform}).
\end{rmk}

\label{rootgroupremark}
The next Lemma follows similarly to Lemma 3.6 of \cite{kumardescent}.

\begin{lem}\label{32}
For $x = (B, \dot{w}_1u_1 B, \dot{w}_2 u_2 B) \in X$, as in equation (\ref{niceform}), 
\[
T_x =  \bigcap_{\alpha \in R_x} \ker(e^\alpha).
\]

\end{lem}

\begin{prf} First, we claim that for any $u \in U^{-1}$
\[
T \cap uBu^{-1} = \cap_{\alpha \in \Rc(u)} \ker(e^\alpha), 
\]
where $e^\alpha: T \rightarrow \C^*$ is the character corresponding to the root $\alpha$. If the claim holds, since for any subset $S \subset R$ and $w \in W$, 
\[
\dot{w}(\bigcap_{\alpha \in S} \ker(e^\alpha))\dot{w}^{-1} = \bigcap_{\beta \in w S} \ker(e^\beta), 
\]
we can conclude that
\[
\begin{array}{ll}
T_x &= T \cap \dot{w_1}u_1 B u_1^{-1} \dot{w_1}^{-1} \cap \dot{w_2}u_2B u_2^{-1} \dot{w_2}^{-1} \\ \\
&= \big[ T \cap \dot{w_1}u_1 B u_1^{-1} \dot{w_1}^{-1}\big]  \cap \big[ T \cap  \dot{w_2}u_2B u_2^{-1} \dot{w_2}^{-1} \big] \\ \\
&= \big[\bigcap_{\alpha \in w_1 \Rc(u_1)} \ker(e^\alpha)\big] \cap \big[\bigcap_{\alpha \in w_2 \Rc(u_2)} \ker(e^\alpha)\big]\\ \\
&= \bigcap_{\alpha \in [w_1 \Rc(u_1) \cup w_2 \Rc(u_2)]} \ker(e^\alpha).
\end{array}
\]

To prove the claim, by Remark \ref{intersectremark} any point in $t \in T \cap uBu^{-1}$ is of the form $t = u t_1 u^{-1}$ for $t_1 \in T$. Since $T$ normalizes $U^{-1}$ and using $B^- = TU^-$ decomposition, we see that $t = t_1$. Thus, $t \in T \cap uBu^{-1}$ if and only if $t = utu^{-1}$, or equivalently $t^{-1}u t = u$. For, $u = u_{\beta_1}(x_1) \cdots u_{\beta_m}(x_m) = \phi(x_1, \dots, x_m)$, this is equivalent to 
\[
 \begin{array}{ll}
  t^{-1}u t& = t^{-1}u_{\beta_1}(x_1) \cdots u_{\beta_m}(x_m)t \\ \\
           & =  t^{-1}u_{\beta_1}(x_1)t t^{-1} \cdots t t^{-1}u_{\beta_m}(x_m)t \\ \\
           & =  u_{\beta_1}(e^{\beta_1}(t^{-1})x_1) \cdots u_{\beta_m}(e^{\beta_m}(t^{-1})x_m) \\ \\
           & = u_{\beta_1}(x_1) \cdots u_{\beta_m}(x_m).
 \end{array}
\]
In other words $\phi(x_1, \dots, x_m) = \phi(e^{\beta_1}(t^{-1})x_1, \dots, e^{\beta_m}(t^{-1})x_m)$. But, as $\phi$ is an isomorphism of varieties, this forces $e^{\beta_i}(t^{-1}) = 1$ for every $x_i \neq 0$. The claim is proved. \end{prf}

\section{Some Conditions for Descent}

We begin this section with a first condition for trivial action of stabilizers.

\begin{lem}\label{26}
Let $x \in X$ be as in equation (\ref{niceform}). If $t \in T_x$, then $t$ acts trivially on $\Lc_x$ if and only if $e^{\lam + w_1 \mu + w_2\nu}(t) = 1$. 
\end{lem}

\begin{prf} Any $t \in T_x = T \cap \dot{w}_1 u_1 B u_1^{-1}\dot{w}_1^{-1} \cap \dot{w}_2 u_2 B u_2^{-1}\dot{w}_2^{-1} $ satisfies equation (\ref{toricform}) for $t_i = \dot{w}_i^{-1}t\dot{w}_i$, $i =1,2$. Also, vectors in $\Lc(\lam, \mu, \nu)_x$ are of the form $[e, z_1] \otimes [ \dot{w}_1 u_1, z_2] \otimes [\dot{w}_2 u_2, z_3]$ for $z_1, z_2, z_3 \in \C^*$.  Thus, the action of $t$ on $\Lc_x$ is  
\[
\begin{array}{l}
t \cdot ([e, z_1] \otimes [ \dot{w}_1u_1, z_2] \otimes [\dot{w}_2u_2, z_3]) \\ \\
= [t, z_1] \otimes [t \dot{w}_1u_1, z_2] \otimes [t\dot{w}_2u_2, z_3] \\ \\
= [t, z_1] \otimes [ (\dot{w}_1u_1t_1 v_1 u_1^{-1} \dot{w}_1^{-1}) \dot{w}_1u_1, z_2] \otimes [(\dot{w}_2 u_2 t_2 v_2 u_2^{-1} \dot{w}_2^{-1})\dot{w}_2u_2, z_3] \\ \\
= [t, z_1] \otimes [ \dot{w}_1u_1t_1 v_1, z_2] \otimes [\dot{w}_2 u_2 t_2 v_2, z_3] \\ \\
= [e, tz_1] \otimes [ \dot{w}_1u_1, t_1 v_1z_2] \otimes [\dot{w}_2 u_2 ,t_2 v_2 z_3] \\ \\
= e^{-\lam}(t)e^{-\mu}(t_1)e^{-\nu}(t_2)[e, z_1] \otimes [ \dot{w}_1u_1, z_2] \otimes [\dot{w}_2u_2, z_3] .
\end{array}
\]
But, by lemma \ref{30} the coefficient is
\[
e^{-\lam}(t)e^{-\mu}(t_1)e^{-\nu}(t_2) = e^{-\lam - w_1\mu -w_2 \nu}(t).
\] \end{prf}

The following lemma appears in \cite{kumardescent} but with the weight lattice $\Lam$ in place of $X(T)$. However, it is easy to observe that the proof given there holds just as well for $X(T)$, and hence we omit a proof. 

\begin{lem}\label{2}
For $S \subset R$ any collection of roots, let $T_S := \cap_{\alpha \in S} \ker(e^\alpha) \subset T$. For any $\lam \in X(T)$, $e^\lam|_{T_S} \equiv 1$ if and only if $\lam \in \Z S$. 
\end{lem}

We state another lemma which is known and hence presented without proof. The proof follows easily by inductively applying simple reflections to fundamental weights. 

\begin{lem}\label{40}
For $\g$ of any type, if $\lam \in \Lam$ and $w \in W$, then $\lam - w \lam \in Q$. 
\end{lem}

Now, define the \textit{lattice associated to $x$} as
\[
 L_x := \Z (R_x).
\]
We now easily derive a necessary condition for descent. 

\begin{prop}\label{100}
If $\Lc$ descends to $X//G$, then $\lam + \mu + \nu \in Q$. 
\end{prop}

\begin{prf} We are always working with the assumption that there is a semistable point. Moreover, we claim the existence of a semistable point of the form $x = (B, \dot{w_1}u_1 B ,  \dot{w_2}u_2 B)$ with $G_x = T_x$. The form of $x$ is a simple matter of selecting a certain point from an orbit. To achieve $G_x = T_x$, we may assume that $G \cdot x$ is closed in $X^{ss}$, for if not, we simply take a point in a minimal dimenion boundary orbit of $x$ within $X^{ss}$. Now apply \cite{dh}, Lemma 3.3.12, to know that $G_x \subset B$ is reductive, and translate by $b \in B$ if necessary according to Proposition \ref{54}. 

Then, $\Lc$ descends only if $T_x$ acts trivially on $\Lc_x$. By Lemma \ref{26}, $T_x$ acts trivially on $\Lc_x$ if and only if $e^{\lam + w_1 \mu + w_2\nu}|_{T_x} =1$, which by Lemma \ref{2} and Lemma \ref{32}, is equivalent to $\lam + w_1 \mu + w_2\nu \in L_x$. But, $L_x \subset Q$ always. So, descent implies $\lam + w_1 \mu + w_2\nu \in Q$ for some $w_1, w_2 \in W$. Now, apply Lemma \ref{40}.
\end{prf}

Define the subset $W^{ss}$ of $W \times W$ as
\[
W^{ss} := \{ (w_1, w_2) \;|\; \exists \; u_1 \in U_{w_1^{-1}}^-\;,\; u_2 \in U_{w_2^{-1}}^-\; \mbox{with}\; x = (B, \dot{w_1}u_1 B ,  \dot{w_2}u_2 B) \in X^{ss}\;\}.
\]
Moreover, define the subset $Z^{ss} \subset (B/B \times G/B \times G/B)\cap X^{ss}$ as
\[
 Z^{ss}:=\{ x = (B, \dot{w_1}u_1 B ,  \dot{w_2}u_2 B) \in X^{ss}\;|\;G_x = T_x\;\}.
\]
Consider a map of sets $\xi:Z^{ss} \rightarrow W^{ss}$, by $\xi(x) = (w_1, w_2)$ when $x = (B, \dot{w_1}u_1 B ,  \dot{w_2}u_2 B)$. We summarize Lemma \ref{32}, Lemma \ref{26}, and Lemma \ref{2} in the following proposition. In this context, if $\xi^{-1}(w_1, w_2) = \emptyset$ for some $(w_1, w_2) \in W^{ss}$, we set $\cap_{x \in \emptyset} L_x = Q$. 

\begin{prop}\label{iffdescentprop}
$\Lc$ descends to $X//G$ if and only if for every $(w_1, w_2) \in W^{ss}$, $\lam + w_1 \mu + w_3 \nu \in \bigcap_{x \in \xi^{-1}(w_1, w_2)} L_x$.  
\end{prop}

We consider the intersection $\Gamma$ of all lattices $ \Z S $ over all subsets $S \subset R^+$ with $\Z S$ of finite index in $Q$. The lattice $\Gamma$ was fully determined in Theorem 3.10 of \cite{kumardescent} by using the Borel-de Siebenthal classification of semisimple subalgebras $\mathfrak{s}$ of $\g$ of maximal rank (see e.g. \cite{wolf}). It is noteworthy to mention that $\Gamma$ is $W$-invariant. We give $\Gamma$ in Table \ref{table:Gammalattice} below.

\begin{table}[h]
\caption{The Lattice $\Gamma$}\label{table:Gammalattice}
\begin{tabular}{| >{$}l<{$}  >{$}l<{$} | }

\hline
 A_\ell \;\; ( \ell \geq 1):& \Gamma = Q. \\
 B_\ell \;\; ( \ell \geq 3):& \Gamma = 2 Q. \\
 C_\ell \;\; ( \ell \geq 2):& \Gamma = 2 \Lam. \\
D_4 \; :&
 \Gamma = \{ n_1\alpha_1 + 2 n_2 \alpha_2 + n_3 \alpha_3 + n_4 \alpha_4\;|\; n_i \in \Z\;, \;\; n_1 + n_3 + n_4 \in 2\Z\}. \\
 D_\ell \;\; ( \ell \geq 5):&
\Gamma = \{ 2n_1\alpha_1 + \cdots + 2 n_{\ell - 2} \alpha_{\ell - 2} + n_{\ell - 1} \alpha_{\ell - 1} + n_\ell \alpha_\ell \;|\; n_i \in \Z \;, \;\; n_{\ell - 1}  + n_\ell \in 2\Z\}. \\
 G_2 \; :& \Gamma = \Z 6 \alpha_1 + \Z 2  \alpha_2 . \\
 F_4 \;:& \Gamma = \Z 6 \alpha_1 + \Z 6  \alpha_2 + \Z 12 \alpha_3 + \Z 12 \alpha_4 . \\
 E_6 \;:& \Gamma = 6 \Lam. \\
E_7 \;:& \Gamma = 12 \Lam. \\
 E_8 \;:& \Gamma = 60 Q . \\\hline
\end{tabular}
\end{table}

First, we establish a non-uniform condition for descent that depends on a lattice that is independent of the GIT class. The proof follows similarly to Theorem 3.9 of \cite{kumardescent} in that we also extend lattices $\Z S$, for $S \subset R^+$, of finite index in $Q$ to lattices of infinite index in $Q$.

\begin{thm}\label{22}
The line bundle $\Lc(\lam, \mu, \nu)$ descends to $X//G$ if $ \lam + w_1\mu + w_2 \nu \in \Gamma$ for every $(w_1, w_2) \in \xi(Z^{ss})$.
\end{thm}

\begin{prf} Assume that $ \lam + w_1\mu + w_2 \nu \in \Gamma$ for every $(w_1, w_2) \in \xi(Z^{ss})$. Let $(w_1, w_2) \in \xi(Z^{ss})$ and $x \in \xi^{-1}((w_1, w_2))$. By Proposition \ref{iffdescentprop}, it suffices to show that $\lam + w_1\mu + w_2 \nu \in L_x$. Whenever $L_x$ is of finite index in $Q$ this certainly holds since $ \Gamma \subset L_x$ by the definition of $\Gamma$. On the other hand, if $L_x$ is infinite index in $Q$, there is a subset $S := \{ \alpha_{i_1}, \cdots , \alpha_{i_k}\}$ of simple roots such that $L_x + \Z S $ is always finite index in $Q$ and $L_x \cap  \Z S = \{0\}$.

Now, some power $\Lc^N$ of $\Lc$ descends to $X//G$ (in fact, for any GIT quotient). Since $G_x = T_x$ and $\Lc^N$ descends to $X//G$, $T_x$ acts trivially on $\Lc^N_x$. Hence, by Lemma \ref{26}, we have $e^{N\lam + Nw_1\mu + Nw_2 \nu}|_{T_x} = 1$, and so $N\lam + Nw_1\mu + Nw_2 \nu \in L_x$ by Lemma \ref{2}. 

Yet, $ L_x + \Z S$ being finite index in $Q$ implies
\[
 \lam + w_1\mu + w_2 \nu \in \Gamma \subset  L_x + \Z S .
\]
Write $\lam + w_1\mu + w_2 \nu  = \alpha_x + \alpha_S$ for $\alpha_x \in L_x$ and $\alpha_S \in \Z S$. Then,
\[
 N(\lam + w_1\mu + w_2 \nu) = N\alpha_x + N\alpha_S \in  L_x .
\]
This implies $ N\alpha_S \in L_x$. But, by choice $L_x \cap  \Z S = \{0\}$. Thus, $\alpha_S = 0$ and we conclude $\lam + w_1\mu + w_2 \nu \in L_x$. \end{prf}

Now, let $\theta$ be the longest root of $\g$ and $d $ be the least common multiple of the coefficients of $\theta$ in terms of the simple roots. For every simple $\g$, both $\theta$ and $d$ are given in Table \ref{table:theta}.

\begin{table}[h]
\centering
\caption{$\theta$ and $d$ for each type of $\g$}\label{table:theta}
\begin{tabular}{| >{$}l<{$}  >{$}l<{$} >{$}l<{$} | }
\hline
 A_\ell: &\theta = \alpha_1 + \alpha_2 + \cdots + \alpha_\ell& d = 1 \\
 B_\ell :& \theta = \alpha_1 + 2\alpha_2 + 2\alpha_3 + \cdots + 2 \alpha_\ell & d = 2 \\
 C_\ell :& \theta = 2\alpha_1 + 2\alpha_2 + \cdots + 2 \alpha_{\ell-1} + \alpha_\ell & d = 2\\
 D_\ell:& \theta = \alpha_1 + 2\alpha_2 + \cdots + 2 \alpha_{\ell-2}+ \alpha_{\ell - 1} + \alpha_\ell & d = 2 \\
 G_2  :&\theta = 3 \alpha_1 + 2 \alpha_2 & d = 6\\
 F_4 :& \theta = 2 \alpha_1 + 3 \alpha_2 + 4 \alpha_3 + 2 \alpha_4 & d = 12 \\
 E_6 :& \theta = \alpha_1 + 2 \alpha_2 + 2\alpha_3 + 3 \alpha_4 + 2 \alpha_5 + \alpha_6 & d =6 \\
E_7 :& \theta = 2\alpha_1 + 2 \alpha_2 + 3\alpha_3 + 4\alpha_4 + 3\alpha_5 + 2 \alpha_6 + \alpha_7& d = 12 \\
 E_8 :&\theta = 2\alpha_1 + 3\alpha_2 + 4\alpha_3 +6\alpha_4 + 5\alpha_5 +4\alpha_6 + 3\alpha_7 + 2\alpha_8& d = 60 \\\hline
\end{tabular}
\end{table}

\begin{lem}\label{almosttherelem}
For $d$ as above, we have $d Q \subset \Gamma$ in all cases. Moreover, when $\lam, \mu, \nu \in  d \Lam$, $\lam + \mu + \nu \in \Gamma$ if and only if $\lam + w_1 \mu + w_2 \nu \in \Gamma$.
\end{lem}

\begin{prf} The first claim is immediate upon inspection. For the second claim, just observe that for $\lam \in  \Lam$ and any $w \in W$, by Lemma \ref{40} we have $d \lam - w d\lam = d(\lam - w \lam) \in d Q  \subset \Gamma$.\end{prf}

If we add an additional assumption on $\lam,\mu,$ and $\nu$, then Lemma \ref{almosttherelem} and Theorem \ref{22} give a sufficient condition for descent that is independent of $W^{ss}$, hence uniform across all GIT classes. 

\begin{thm}\label{56}
For any $\lam, \mu, \nu  \in X(T)^{++}$, if $\lam, \mu, \nu \in  d \Lam$ and $\lam + \mu +  \nu \in \Gamma$, then $\Lc(\lam, \mu, \nu)$ descends to $X//G$. 
\end{thm}

\begin{rmk}
Although we technically only require that $ \mu, \nu \in  d \Lam$ and $\lam + \mu +  \nu \in \Gamma$. The following result shows that this condition as stated is not much of a loss of generality. 
\end{rmk}

\begin{lem}
In all cases except $G_2 $ and $F_4$, $\Gamma \subset d \Lam $. If $G$ is of type $G_2$ or $F_4$, then, $d \Lam \subset \Gamma$.
\end{lem}

\begin{prf} The statement is clear in all cases except for $G_2$, $F_4$, and $D_\ell$, $\ell \geq 4$. For $G_2$ and $F_4$, simply note that $\Lam = Q$, and since $d Q \subset \Gamma$, we are done. For $D_\ell$, clearly $2Q \subset 2 \Lam$. When $\ell \geq 5$, it suffices to show $m \alpha_{\ell-1} + n \alpha_\ell  \in 2 \Lam $ if $n + m $ is even. Since, $2 \varpi_{\ell-1} - 2\varpi_\ell = \alpha_{\ell-1} - \alpha_\ell$ and also $2 \alpha_{\ell -1} \in 2 \Lam$, we have
\[
m \alpha_{\ell-1} + n \alpha_\ell = (m+n)\alpha_{\ell-1} - n(\alpha_{\ell-1} - \alpha_\ell) \in 2 \Lam
\]
whenever $m + n $ is even. When $\ell = 4$, it suffices to show $a\alpha_1 + b \alpha_3 + c \alpha_4 \in 2 \Lam$ if $a + b + c $ is even. But, $2 \varpi_1 - 2 \varpi_3 = \alpha_1 - \alpha_3 \in 2 \Lam$. Hence
\[
(a + b + c) \alpha_1 - (b+ c) (\alpha_1 - \alpha_3) - c(\alpha_3 - \alpha_4) 
=a\alpha_1 + b \alpha_3 + c \alpha_4 \in 2 \Lam
\]
whenever $a + b + c $ is even. \end{prf} 

Since $\Gamma$ is $W$-invariant, $\lam - w \lam \in \Gamma$ for any $\lam \in \Gamma$ and any $w \in W$. From this fact, we get a slight strengthening of Theorem \ref{56} for $G_2$ and $F_4$. Of course, the following corollary holds in all cases, yet the corollary is weaker than Theorem \ref{56} except for $G_2$ and $F_4$.

\begin{corl}\label{nonmaindescent}
If $\lam, \mu, \nu \in \Gamma$ then $\Lc(\lam, \mu, \nu)$ descends to $X//G$. 
\end{corl}

\begin{rmk} We see that descent is equivalent to a list of characters lying in certain sublattices of the root lattice by Proposition \ref{iffdescentprop}. Yet, it is poorly understood which characters are relevant (i.e., determining $W^{ss}$), and which sublattice the relevant characters must lie in (i.e., which $\Rc(u_1)$ and $\Rc(u_2)$ can occur with some $(w_1, w_2) \in W^{ss}$). 

To add difficulty to the issue, it is unclear whether or not this gives a condition that is uniform accross all GIT classes; or, if some GIT classes have a descent condition that is unrelated to the other GIT classes. We briefly expand upon this. Recall that the collection of ample $G$-equivariant line bundles on $X$ is equipped with an equivalence relation where line bundles are related if both their semistable and stable loci agree. The equivalence classes are called GIT classes. GIT classes form convex cones, and the GIT classes of maximal dimension are those with empty stable locus (see \cite{dh} for a thorough treatment). Since $X^{ss}$ is constant within a GIT class, it is clear that $W^{ss}$ and $Z^{ss}$ are also constant within a GIT class. Thus, Proposition \ref{iffdescentprop} gives a uniform descent condition accross a single GIT class. However, since $X^{ss}$ changes as we change GIT classes, it is not clear how $W^{ss}$ and $Z^{ss}$ will change and how this will affect the corresponding descent condition.

In section \ref{counterexample} we see an example of this issue. The triple $(2\rho, \rho, \rho)$ has $\{(w_0, w_0)\} = W^{ss}$, where $w_0$ is the longest Weyl group element. However, since $(B, w_0B, w_0B) \in X^{ss}$, we have $\bigcap_{x \in \xi^{-1}(w_0, w_0)} L_x = \{0\}$. One can check using the inequalities of \cite{belkalekum} that any GIT class in the Cartan component (see \cite{kumartensor} for definition), and hence $(2\rho, \rho, \rho)$, lives in the boundary of $C^G(X)$. On the other hand, if $\Lc$ is in a GIT class of maximal dimension, and so also in the interior of $C^G(X)$, then one can show that $L_x$, for any $x \in \xi^{-1}(w_1, w_2)$, is of finite index in $Q$. Despite this apparently vast difference, the proof of Theorem \ref{22} shows that any $L_x$ of infinite index in $Q$ can be `replaced' by a lattice of finite index in $Q$. So, it remains unclear how different descent conditions truly are as we vary the GIT classes of $\Lc$. 

Regardless, by allowing the additional assumptions of Theorem \ref{56}, we obtain sufficient conditions for descent that are uniform accross all GIT classes. Although it is clear from section \ref{counterexample} that our sufficient conditions for descent are not optimal for some GIT class, it is not clear if our sufficient conditions are the optimal conditions uniform accross all GIT classes. 

\end{rmk}

\section{A Few Short Examples}

\subsection{Explicit computation for $SL(2)$ and $SL(3)$}

Here, we compute sufficient conditions for descent by hand for $G= SL(2)$ and $G = SL(3)$, and observe that these conditions match those given in Theorem \ref{56}. Increasing the rank of these examples increases the difficulty of `by hand' computations significantly, so these examples are restricted to very low rank cases. The following lemma follows as in lemma 3.5 of \cite{kumardescent}, and is helpful in our examples. 

\begin{lem}\label{33}
Let $x \in X^{ss}$ be as in equation (\ref{niceform}). When $T_x$ is a divisible group, $e^{\lam + w_1 \mu + w_2\nu}\mid_{T_x} \equiv 1$.  
\end{lem}

\begin{prf} For $x \in X^{ss}$, there exists $\sigma \in H^0(X, \Lc^N)^G$, for some $N> 0$, such that $\sigma (x) \neq 0$. Recall that $\Lc^N = \Lc(N\lam) \boxtimes \Lc(N\mu) \boxtimes \Lc(N\nu)$ and let 
\[
\sigma(x) = [e, z_1] \otimes [ \dot{w}_1u_1, z_2] \otimes [\dot{w}_2u_2, z_3].
\]
$G$-invariance of $\sigma$ implies $  \sigma (x) = (t \cdot \sigma)(x) = t\cdot \sigma(t^{-1} \cdot x) = t\cdot \sigma( x)$ when $t \in T_x$. Then, by Lemma \ref{26}
\[
t\cdot \sigma( x) = e^{-N\lam - w_1N\mu -w_2N \nu}(t)[e, z_1] \otimes [ \dot{w}_1u_1, z_2] \otimes [\dot{w}_2u_2, z_3].
\]
Hence 
\[
e^{-N\lam - Nw_1\mu -Nw_2 \nu}|_{T_x} \equiv 1.
\]
Then, for any $t \in T_x$, divisibility implies $s^N = t$ for some $s \in T_x$. So, 
\[
e^{-\lam - w_1\mu -w_2 \nu}(t) = e^{-\lam - w_1\mu -w_2 \nu}(s^N)=e^{-N\lam - Nw_1\mu -Nw_2 \nu}(s) =1.
\]\hfill\end{prf}

For $SL(2)$, it is easy to see that the different possibilities for $T_x$ are $ T $ and $\{\pm I\}$. Hence, the only case where $T_x$ is not a divisible group is $\{\pm I\}$. If, 
\[
\lam = b_1 \varpi_1, \;\; \mu = b_2 \varpi_1, \;\; \nu = b_3 \varpi_1,
\]
then, for any $w_1, w_2 \in W$
\[
e^{-\lam - w_1 \mu - w_2\nu}(-I) = (-1)^{b_1 + b_2 + b_3}.
\]
Hence, we have descent if $2 \mid b_1 + b_2 + b_3$, which is equivalent to $\lam + \mu + \nu \in Q$.

For the $SL(3)$ case, a more involved computation shows that the only possibilities for $T_x$ are $T$, $\C^*$, and the three element group generated by 
\[
\zeta I = \begin{pmatrix}
\zeta & 0 & 0 \\
0 & \zeta & 0 \\
0 & 0 & \zeta
\end{pmatrix},
\]
where $\zeta $ is a primitive cube root of unity. The first two cases are divisible. 
If, 
\[
\lam = ( a_1- b_1)\varpi_1 + b_1\varpi_2, \;\; \mu =( a_2- b_2)\varpi_1 + b_2\varpi_2, \;\; \nu =( a_3- b_3)\varpi_1 + b_3\varpi_2,
\]
then, for any $w_1, w_2 \in W$,
\[
e^{-\lam - w_1 \mu - w_2\nu}(\zeta I) = \zeta^{\sum_{i=1}^3 a_i + b_i}.
\] 

In particular, this is trivial if $3 | \sum_{i=1}^3 a_i + b_i$, which is again equivalent to $\lam + \mu + \nu  \in Q$.

\subsection{A Counter Example to ``Sufficient is Neccesary"}\label{counterexample}
Here we consider the triple $(\lam, \mu, \nu) = (2\rho, \rho, \rho)$ for $\rho = \sum_{\alpha \in R^+} \frac{1}{2} \alpha$. We show that $\Lc(2\rho, \rho, \rho)$ always descends to $X//G$. Except for $\g  = \mathfrak{sl}_\ell$, $(2\rho, \rho, \rho)$ does not satisfy one or both (depending on $\g$) of the conditions given in Theorem \ref{56}. One can make similar arguments for any triple $(-w_0(\lam + \mu),\lam, \mu)$, where $w_0$ is the longest element of $W$, i.e., the so-called Cartan components (see e.g., \cite{kumartensor}). However, we consider only $(2\rho, \rho, \rho)$ for simplicity of exposition. 

\begin{prop}
 The line bundle $\Lc(2\rho, \rho, \rho)$ always descends to $X//G$. Moreover, in this case, $X^{ss} = G \cdot (B, w_0U^-B, w_0U^-B)$.
\end{prop}

\begin{prf} Recall that $V(N\rho)^* \cong V( N\rho)$. Let $\phi_{N\rho} \in V(N \rho)^*$ be a highest weight vector dual to a lowest weight vector $v_{ - N \rho} \in V(N \rho)$. Then, we have the canonical isomorphism $V(N\rho)^* \cong V(- w_0 N\rho) = V(N \rho)$ obtained by extending $G$-linearly $\phi_{N\rho} \mapsto v_{N \rho}$, where $v_{N \rho}$ is a highest weight vector of $V(N \rho)$, and noting that $-w_0 \rho = \rho$. Then, by \cite{kumartensor}, Lemma 3.1, we know 
\[
 \begin{array}{ll}
  \Hom_G(V(2N\rho), V(N\rho) \otimes V(N\rho)) &\cong [V(2N\rho)^* \otimes V(N\rho) \otimes V(N\rho)]^G\\\\
  &\cong [V(2N\rho)^* \otimes V(N\rho)^* \otimes V(N\rho)^*]^G\\\\
    &\cong H^0(X, \Lc(2\rho, \rho, \rho)^N)^G\\\\
 \end{array}
\]
is one dimensional. It suffices to determine a single non-zero section $\sigma_N \in   H^0(X, \Lc(2\rho, \rho, \rho)^N)^G$ in order to determine $X^{ss}$. We construct such $\sigma_N$. 

Let $\psi^\circ$ be the equivariant embedding of $V(2N\rho) \hookrightarrow V(N\rho)\otimes V(N\rho)$ given by the Cartan component as mentioned above. So, $\psi^\circ$ extends $G$-linearly $\psi^\circ(v_{2N\rho}) = v_{N\rho} \otimes v_{N\rho}$. Composing $\psi^\circ$ with the isomorphisms $V(N\rho) \cong V(N\rho)^*$ above gives $\psi \in \Hom_G(V(2N\rho), V(N\rho)^* \otimes V(N\rho)^*)$. So, $\psi$ extends $G$-linearly $\psi(v_{2N\rho}) = \phi_{N\rho} \otimes \phi_{N\rho}$ . 

Now, let $\{w_{\gamma}^i\}$, $\{v_\mu^i\}$ be bases for $V(2N \rho)$ and $V(N\rho)$, respectively, where the basis vector $v_\mu^i$ is taken in weight space $\mu$, with $i$ indexing basis vectors within each weight space, and similarly for the $w_\gamma^i$. Moreover, let $\{(w_\gamma^i)^*\}$, $\{(v_{\mu}^i)^*\}$, be respective dual bases for $V(2N\rho)^*$ and $V(N\rho)^*$ where the weight space subscript indicates the weight space the vector is dual to. For example, $(v_{-N\rho})^* = \phi_{N\rho}$. Then, applying the usual isomorphism $  \Hom_G(V,W) \cong (V^* \otimes W)^G$ to $\Hom_G(V(2N\rho), V(N\rho)^* \otimes V(N\rho)^*)$, gives
\begin{equation}\label{sectioneq}
\begin{array}{ll}
\psi &\mapsto \sum_{\gamma,i} (w_\gamma^i)^* \otimes \psi(w_\gamma^i) \\\\
&= (w_{2N \rho})^* \otimes \psi(w_{2N \rho}) + \sum_{\gamma < 2N \rho,i} (w_\gamma^i)^* \otimes \psi(w_\gamma^i) \\\\
&= (w_{2N \rho})^* \otimes \phi_{N\rho} \otimes \phi_{N\rho} + \sum_{\gamma < 2N \rho,i} (w_\gamma^i)^* \otimes \psi(w_\gamma^i).
\end{array} 
\end{equation}
This gives a non-zero vector in $[V(2N\rho)^* \otimes V(N\rho)^* \otimes V(N\rho)^*]^G$.

Recall, for a simple tensor $\varphi_1\otimes \varphi_2 \otimes \varphi_3 \in V(\lam)^* \otimes V(\mu)^* \otimes V(\nu)^*$, the Borel-Weil isomorphism, which we denote $\beta$, gives the section 
\[
 \beta( \varphi_1\otimes \varphi_2 \otimes \varphi_3)(g_1B, g_2B, g_3B) = [g_1, \varphi_1(g_1 v_\lam)]\otimes [g_2, \varphi_2(g_2 v_\mu)]\otimes [g_3, \varphi_3(g_3 v_\nu)],
\]
where $v_\lam, v_\mu$, and $v_\mu$ are highest weight vectors of the respective representations. Applying the Borel-Weil isomorphism to equation (\ref{sectioneq}), we get a $G$-invariant section $\sigma_N$. We write $\sigma_N= \sigma_{\mbox{\tiny top}} + \sigma_{\mbox{\tiny low}}$, where $\sigma_{\mbox{\tiny top}}$ is the image of $(w_{2N \rho})^* \otimes \phi_{N\rho} \otimes \phi_{N\rho}$, and $\sigma_{\mbox{\tiny low}}$ is the image of $\sum_{\gamma < 2N \rho,i} (w_\gamma^i)^* \otimes \psi(w_\gamma^i)$. 

Since semistability depends only on $G$-orbits, to determine the zero set of $\sigma_N$ it suffices to consider only points of the form $(B, g_1 B, g_2B)$. Since, $( w_\gamma^i)^*( v_{2N\rho})$ is non-zero if and only if $\gamma = 2N\rho$, for all $\gamma < 2N \rho$,
\[
\beta((w_\gamma^i)^* \otimes \psi(w_\gamma^i))(B, g_1 B, g_2B) = K \cdot  (w_\gamma^i)^*( v_{2N \rho}) = 0, 
\]
where $K \in \C$. Hence, $\sigma_{\mbox{\tiny low}}(B, g_1 B, g_2B) = 0$. 

Then,
\[
\sigma_N(B, g_1 B, g_2B)= \sigma_{\mbox{\tiny top}}(B, g_1 B, g_2B) = [e,1] \otimes [g_1, \phi_{N\rho}(g_1 v_{N\rho}^+)]  \otimes [g_2, \phi_{N\rho}(g_2v_{N\rho}^+)].
\]
This gives the following description.
\[
X^{ss} = G \cdot \{(B, g_1B, g_2B)\;|\; [g_1 v_{N\rho}^+]_{-N\rho} \neq 0,\;\;\mbox{and}\;\;[g_3 v_{N\rho}^+]_{-N\rho} \neq 0\;\;\;\mbox{for some}\;\; N > 0\},
\]
where $[w]_\mu$ denotes the $\mu$-weight space component of any vector $w$. In fact, using Bruhat decomposition one can readily see that  see that this set corresponds to  
\[
 G \cdot (B, Bw_0B, Bw_0B) = G \cdot (B, w_0U^-B, w_0U^-B). 
\]

So, we see that $W^{ss}= \{(w_0, w_0)\}$, and so descent easily follows from Lemma \ref{26} since $e^{2\rho + w_0 \rho + w_0 \rho} = e^0 = 1$. 
\end{prf}

\bibliographystyle{abbrv}
\bibliography{descentbiblio}

\end{document}